\documentclass[12pt, reqno]{amsart}
\usepackage{amssymb,a4}
\usepackage{amsmath,amsthm,amsfonts,amssymb}

\def\a{\alpha}
\def\b{\beta}

\def\D{\Delta}
\def\UU{{\mathcal U}}

\def\LL{{\mathcal L}}
\def\Der{{\rm Der}}
\def\VV{\mathcal {V}}
\def\Inn{{\rm Inn}}
\def\Ker{{\rm Ker}}
\def\Im{{\rm Im}}
\def\v{\varphi}
\def\ssc{\scriptscriptstyle}

\def\ma{\mathbb}
\def\rar{\rightarrow}
\def\bs{\backslash}
\def\vs{\vspace*}
\def\even{{\bar0}}
\def\odd{{\bar1}}

\def\ni{\noindent}
\def\N{\mathbb{N}{\ssc\,}}
\def\Z{\mathbb{Z}{\ssc\,}}
\def\C{\mathbb{C}{\ssc\,}}

\def\QED{\hfill$\Box$}

\numberwithin{equation}{section}
\newtheorem{theo}{Theorem}[section]
\newtheorem{defi}[theo]{Definition}
\newtheorem{rema}[theo]{Remark}
\newtheorem{lemm}[theo]{Lemma}
\newtheorem{prop}[theo]{Proposition}
\newtheorem{clai}{Claim}

\title{Lie superbialgebra structures on the $N=2$ superconformal Neveu-Schwarz algebra}

\author{Dong Liu}
\address{Department of Mathematics, Huzhou Teachers College, Zhejiang Huzhou, 313000, China}
\email{liudong@hutc.zj.cn}

\author{Liangyun Chen}
\address{Department of Mathematics, Northeast Normal University,
Changchun, 130024, China} \email{Corresponding
author:chenly640@nenu.edu.cn}

\author{Linsheng Zhu}
\address{Department of Mathematics, Changshu Institute of
Technology, Jiangsu Changshu, 215500, China}
\email{lszhu@cslg.edu.cn}

\date{}
\begin{document}
\maketitle

{\small
\parskip .005 truein
\baselineskip 3pt \lineskip 3pt

\noindent{{\bf Abstract.} In this paper, Lie superbialgebra
structures on the $N=2$ superconformal Neveu-Schwarz algebra are
considered by a very simple method. We prove that every Lie superbialgebra structure on the
algebra is triangular coboundary. \vs{5pt}

}

}
\section{Introduction}
\setcounter{section}{1}\setcounter{equation}{1}

The superconformal algebras, closely related the conformal field
theory and the string theory, play important roles in both
mathematics and physics supplying the underlying symmetries of
string theory. It is well known that the $N=2$ superconformal
algebras include four sectors: the Neveu-Schwarz sector, the Ramond
sector, the topological sector and the twisted sector. There are
many researches about these algebras (see \cite{CK,VK,DB,EG,KL,SS}
and the reference cited therein). All sectors are closely related to
the Virasoro algebra and the super-Virasoro algebra which play great
roles in the two-dimensional conformal filed.

Lie bialgebras was introduced in 1983 by Drinfeld (see \cite{D1,D2})
during the process of investigating quantum groups. There appeared
several papers on Lie bialgebras and Lie superbialgebras (e.g.,
\cite{M,N,NT,Y}). In \cite{M}--\cite{NT}, the Lie bialgebra
structures on Witt and Virasoro algebras were investigated, which
are shown to be triangular coboundary. Moreover, the Lie bialgebra
structures on the one-sided Witt algebra were completely classified.

In \cite{FL, LFZ,YS}, the Lie superbialgebra structures on the $N=2$
superconformal twisted, topological, Ramond algebras were
investigated case by case with complicated computations.
Clearly, the above superalgebras and the $N=2$ superconformal Neveu-Schwarz algebra all include the centerless twisted  Heisenberg-Virasoro algebra (see Section 2.2 for the definition) as a subalgebra. It is naturally to consider such questions based on the related results of the twisted Heisenberg-Virasoro algebra.

As the above considerations, we study the Lie superbialgebra
structures on the $N=2$ superconformal Neveu-Schwarz algebra $\mathcal L$ in this paper. Successfully, we
provide a very simple method (no need some complicated calculations as in \cite{FL, LFZ,YS}) to
determine the structure of $H^1(\mathcal L, \mathcal L\otimes \mathcal L)$ based on
such results for the twisted Heisenberg-Virasoro algebra in
\cite{LP}. With our considerations, the original proofs in \cite{FL,
LFZ,YS} can be great simplified. Throughout the paper, we denote by
$\Z$, $\C$ the set of all integers, complex numbers, $\Z_+$ (resp
$\Z^*$) the set of all nonnegative (resp. nonzero) integers.

\section{Preliminaries}
\setcounter{section}{2}\setcounter{equation}{0}

\subsection{Lie super-bialgebras}

Firstly, let us recall some related definitions. Let
$L= L_\even\oplus  L_\odd$ be a vector space over $\C$. If $x\in L_{[x]}$, then we say that $x$ is
homogeneous of degree $[x]$ and we write ${\rm deg}x=[x]$. Denote by
$\tau$ the {\it super-twist map} of $ L \otimes  L$, i.e.,
\begin{eqnarray*}
\tau(x\otimes y)= (-1)^{[x][y]}y\otimes x,\ \ \,\forall\,\,x,y\in
 L.
\end{eqnarray*}
For any $n\in\N$, denote by $ L^{\otimes n}$ the tensor product of
$n$ copies of $ L$ and $\xi$ the {\it super-cyclic map} cyclically
permuting the coordinates of $ L^{\otimes3}$, i.e.,
\begin{eqnarray*}
&&\xi=(1\otimes\tau)\cdot(\tau\otimes1): \,x_{1} \otimes
x_{2} \otimes x_{3}\mapsto (-1)^{[x_1]([x_2]+[x_3])}x_{2} \otimes
x_{3} \otimes x_{1},\ \ \forall\,\,x_i\in  L,\,i=1,2,3,
\end{eqnarray*}
where $1$ is the identity map of $ L$. Then the definition of a
Lie superalgebra can be described in the following way: A {\it Lie
superalgebra} is a pair $( L,\v)$ consisting of a vector space
$ L= L_\even\oplus  L_\odd$ and a bilinear map $\v : L
\otimes L \to L$ satisfying:
\begin{eqnarray}\label{cond1}\begin{array}{lll}
&&\v( L_{\bar{i}}, L_{\bar{j}})\subset  L_{\bar{i}+\bar{j}},\vs{4pt}\\
&&\Ker(1 -\tau) \subset \Ker\,\v,\vs{4pt}\\
&&\v \cdot (1  \otimes \v ) \cdot (1 +
\xi +\xi^{2})=0.
\end{array}\end{eqnarray}
Meanwhile, the definition of a Lie super-coalgebra can be described
in the following way: A {\it Lie super-coalgebra} is a pair
$( L,\D)$ consisting of a vector space $ L= L_\even\oplus
 L_\odd$ and a linear map $\D:  L \rar  L \otimes  L$
satisfying:
\begin{eqnarray}\label{cond2}\begin{array}{lll}
&&\D( L_{\bar{i}})\subset
\mbox{$\sum\limits_{\bar{j}\in{\Z_2}}$} L_{\bar{j}}\otimes
 L_{\bar{i}-\bar{j}},\vs{4pt}\\
&& \Im\,\D \subset \Im(1 - \tau),\vs{4pt}\\
&&(1 \otimes1 + \xi +\xi^{2}) \cdot (1 \otimes
\D) \cdot \D =0.
\end{array}\end{eqnarray}
Now one can give the definition of a Lie superbialgebra, which is a
triple $( L, \v, \D )$ satisfying:
\begin{eqnarray*}
&&\mbox{(i)}\ \ ( L, \v){\mbox{ is a Lie superalgebra}},\\[2pt]
&&\mbox{(ii)}\ \ ( L, \D){\mbox{ is a Lie super-coalgebra}},\\[2pt]
&&\mbox{(iii)}\ \ \D  \v (x\otimes y)=x \ast  \D y - (-1)^{[x][y]}y
\ast  \D x\ \ \forall\,\,x, y \in  L,
\end{eqnarray*} where
the symbol ``$\ast$'' means the {\it adjoint diagonal action}
\begin{eqnarray}\label{e-diag}
x\ast  (\mbox{$\sum\limits_{i}$}{a_{i} \otimes b_{i}}) =
\mbox{$\sum\limits_{i}$} ( {[x, a_{i}] \otimes b_{i} +
(-1)^{[x][a_i]}a_{i} \otimes [x, b_{i}]}),\ \ \forall\,\,x, a_{i},
b_{i} \in  L,
\end{eqnarray}
and in general $[x,y]=\v(x\otimes y)$ for $x,y \in  L$.

Denote by $\UU( L )$ the {\it universal enveloping algebra} of
$ L$ and $A\bs B=\{x\,|\,x\in A,x\notin B\}$ for any two sets $A$
and $B$. If $r =\mbox{$\sum\limits_{i}$}{a_{i} \otimes b_{i}} \in
 L \otimes  L $, then the following elements are in
$\UU( L)\otimes \UU( L)\otimes \UU( L)$
\begin{eqnarray*}
&&r^{12} =\mbox{$\sum\limits_{i}$}{a_{i} \otimes b_{i} \otimes 1}=r\otimes1,\ r^{23} =\mbox{$\sum\limits_{i}$}{1 \otimes a_{i}
\otimes
b_{i}}=1\otimes r,\\
&&r^{13}=\mbox{$\sum\limits_{i}$} {a_{i} \otimes 1 \otimes b_{i}}
=(1\otimes\tau)(r\otimes 1)=(\tau\otimes1)(1\otimes r),
\end{eqnarray*}
while the following elements are in $ L\otimes L\otimes L$
\begin{eqnarray*}
&&[r^{12} ,r^{23}]=\mbox{$\sum\limits_{i,j}$}a_{i} \otimes[
b_{i},a_{j}] \otimes
b_{j},\\
&&[r^{12} ,r^{13}]=\mbox{$\sum\limits_{i,
j}$}(-1)^{[a_j][b_i]}[a_{i}, a_{j}] \otimes
b_{i} \otimes b_{j},\\
&&[r^{13} ,r^{23}]=\mbox{$\sum\limits_{i,j} (-1)^{[a_j][b_i]}$}a_{i}
\otimes a_{j} \otimes[ b_{i}, b_{j}].
\end{eqnarray*}

\begin{defi} (i) A {\it coboundary superbialgebra} is a quadruple $( L , \v,
\D,r),$ where $( L , \v, \D)$ is a Lie superbialgebra and $r \in
\Im(1  - \tau) \subset  L \otimes  L $ such that $\D=\D_r$ is a {\it
coboundary of $r$}, i.e.,
\begin{eqnarray}\label{e-D-r}
\D_r(x)=(-1)^{[r][x]}x\ast r,\ \ \forall\,\,x \in  L.
\end{eqnarray}
(ii)\ \ A coboundary Lie superbialgebra $( L , \v,\D, r)$ is called
{\it triangular} if it satisfies the following {\it classical
Yang-Baxter Equation}
\begin{eqnarray}\label{e-CYBE}
c(r):=[r^{12}, r^{13}] +[r^{12} , r^{23}] +[r^{13} , r^{23}]=0.\ \ \
{\rm (CYBE)}
\end{eqnarray}
\end{defi}

Let $V=V_{\bar0}\oplus V_{\bar1}$ be an $ L$-module where
$ L= L_{\bar0}\oplus L_{\bar1}$. A $\Z_2$-homogenous linear map
$d: L\to V$ is called a {\it homogenous derivation of degree
$[d]\in\Z_2$}, if $d( L_i)\subset V_{i+[d]}\
\,(\forall\,\,i\in\Z_2)$,
\begin{eqnarray}\label{e-der}
&&d([x,y])=(-1)^{[d][x]}x\ast d(y)-(-1)^{[y]([d]+[x])}y\ast  d(x),\
\ \forall\,\,x,y\in  L.
\end{eqnarray}
Denote by $\Der_{\bar{i}}( L,V)\ \,(\,i=0,1)$ the set of all
homogenous derivations of degree $\bar{i}$. Then the set of all
derivations from $ L$ to $V$
$\Der( L,V)=\Der_{\bar0}( L,V)\oplus\Der_{\bar1}( L,V)$. Denote
by $\Inn_{\bar{i}}( L,V)\ \,(\,i=0,1)$ the set of {\it homogenous
inner derivations of degree $\bar{i}$}, consisting of $a_{\rm inn},$
$a\in V_{\bar{i}}$, defined by
\begin{equation}
\label{e-inn} a_{\rm inn}:x\mapsto (-1)^{[a][x]}x\ast a,\ \
\forall\,\,x\in  L.
\end{equation}
Then the set of inner derivations
$\Inn( L,V)=\Inn_{\bar0}( L,V)\oplus\Inn_{\bar1}( L,V)$.

Denote by $H^1( L,V)$ the {\it first cohomology group} of $ L$
with coefficients in $V$. Then
\begin{equation*}
H^1( L,V)\cong\Der( L,V)/\Inn( L,V).
\end{equation*}
An element $r$ in a superalgebra $ L$ is said to satisfy the {\it
modified Yang-Baxter equation} if
\begin{equation}\label{e-MYBE}
x\ast  c(r)=0,\ \ \forall\,\,x\in  L.\ \ \ {\rm (MYBE)}
\end{equation}

The following result for the non-super case can be found in
\cite{NT} while its super case can be found in \cite{YS}.
\begin{lemm}\label{0502m01}
Let $L$ be a Lie superalgebra,
$r\in\Im(1-\tau)\subset L\otimes L$ with $[r]=\bar0$.
Then
\begin{equation}\label{e-p2.1}
(1+ \xi + \xi^{2}) \cdot (1\otimes \D_r) \cdot \D_r (x) = x
\ast c (r),\ \ \forall\,\,x \in  L.
\end{equation}
Thus $(L,[\cdot,\cdot],\D_r)$ is a Lie superbialgebra if and only if
$r$ satisfies $(MYBE)$ $({\rm see} (\ref{e-MYBE}))$.
\end{lemm}

\subsection{The $N=2$ superconformal Neveu-Schwarz algebra}

As a vector space over $\C$, the $N=2$ superconformal Neveu-Schwarz algebra $\hat{\mathcal L}$ has
a basis $\{L_m, I_m, G_r^\pm\mid m\in\Z, r\in\Z+1/2\}$, with the following relations:
\begin{eqnarray}\label{LieB}
&& [L_m,L_n]=(m-n)L_{n+m}+{1\over12}(m^3-m)C,\vs{6pt}\\
&&[I_m,I_n]={1\over3}m\delta_{m+n,0}C,\ \ \ \ \ \ \ \ \ \ \ \ \ \ \ \ [L_m, I_n]=-nI_{m+n},\vs{6pt}\\
&&
[L_m,G_r^\pm]=(\frac{m}{2}-r)G_{r+m}^\pm,\ \ \ \ \ \ \ \ [I_m,G_r^\pm]=\pm G_{m+r}^\pm,\vs{6pt}\\
&&[G_r^+,G_s^-]=2L(r+s)+(r-s)I(r+s)+{1\over3}(r^2-{1\over4})\delta_{r+s,0}C,\\
&&[G_r^\pm,G_s^\pm]=0,
\end{eqnarray}
for $m,n\in\Z,\,r,s\in\frac{1}{2}+\Z$.

Denote by ${\mathcal L}$ the centerless Lie superalgebra of $\hat{\mathcal L}$, then $\hat{\mathcal L}$ is the
universal central extension of ${\mathcal L}$.
Obviously, $\LL$ is ${\mathbb{Z}}_2$-graded: ${\mathcal L}={\mathcal
L}_{\overline{0}}\oplus{\mathcal L}_{\overline{1}},$ with
\begin{eqnarray}\label{0426a002}
{\mathcal L}_{\bar{0}}=\mbox{span}_{\ma C}\{
L_m,\,I_m|\,m\in\Z\},\ \ \ \ {\mathcal
L}_{\bar{1}}=\mbox{span}_{\ma C}\{G_r^\pm\,|\,r\in\Z+\frac{1}{2}\}.
\end{eqnarray}
Clearly, ${\mathcal L}_{\bar{0}}$ is the centerless twisted
Heisenberg-Virasoro algebra, and $W=\C\{L_m\mid m\in\Z\}$ is the
Witt algebra. Moreover, ${G}^\pm=\C\{G_{r}^\pm\mid r\in\Z+\frac12\}$
are ${\mathcal L}_{\bar0}$-modules.

For any $\a,\a^\dag,\b,\b^\dag,\gamma,\gamma^\dag\in\C$, one can easily
verify that the linear map $\varrho:{\mathcal L}_{\bar{0}}\to{\mathcal L}_{\bar{0}}\otimes {\mathcal L}_{\bar{0}}$ defined below is a
derivation:
\begin{eqnarray}\varrho(L_n)&=&(n\a+\gamma)I_0\otimes
I_n+(n\a^\dag+\gamma^\dag)I_n\otimes I_0,\nonumber\\
\varrho(I_n)&=&\b I_0\otimes I_n+\b^\dag I_n\otimes I_0,
\ n\in \Z.\nonumber
\end{eqnarray}
Denote $\mathcal D$ the vector space spanned by the such elements
$\varrho$ over $\C$. From Theorem 3.2 and Corollary 4.5 in
\cite{LP}, we have following propositions.

\begin{prop}\label{prop1}\cite{LP}
$H^1({\mathcal L}_{\bar{0}},{\mathcal L}_{\bar{0}}\otimes {\mathcal L}_{\bar{0}})=\mathcal D$.
\end{prop}

\begin{prop}\label{der1}\cite{LP}
$H^1({\mathcal L}_{\bar{0}}, {G}^\pm\otimes {G}^\pm)=H^1({\mathcal L}_{\bar{0}}, {G}^\pm\otimes {G}^\mp)=0$.
\end{prop}

\newpage
The following lemma can be obtained by using the similar techniques
of \cite[Lemma 2.2]{SSu}.
\begin{lemm} \label{lemma2}
Regarding $\LL^{\otimes n}$ as an $\LL$-module under the adjoint
diagonal action of $ \LL $, if $r\in \LL^{\otimes n}$ such that $x\ast
r=0,\ \forall\,\,x\in \LL$, then $r=0$.
\end{lemm}

\noindent{\bf Proof.} It is easy to see that $\mathcal L^{\otimes
n}$ is $\frac{1}{2}\Z$-graded by
\par
$$\mathcal L^{\otimes
n}_{p}=\mbox{$\sum\limits_{p_1+p_2+\cdots+p_n=p}\mathcal L_{p_1}\otimes
\mathcal L_{p_2}$}\otimes\cdots \otimes \mathcal L_{p_n}, \ \ \forall \ \ p, \
p_i\in \frac{1}{2}\Z, \ \ i=1,2,\cdots,n .
$$Write $r=\sum_{p\in\frac{1}{2}\Z}r_{p}$ as a finite sum with $r_p\in
\mathcal L^{\otimes n}_p$. By hypothesis, $L_0\ast r=0$, which
implies $r\in \LL_0$.

So
\begin{equation}r=\sum\limits_{r_1+r_2+\cdots+r_n=0}c_{r_1,r_2,\cdots
r_n}E_{r_1}\otimes E_{r_2}\otimes\cdots\otimes
E_{r_n}\label{equ-r}\end{equation} for some $c_{r_1,r_2,\cdots
r_n}\in \C$ and $E_{r_i}\in\mathcal L_{r_i}$ with
$r_i\in\frac{1}{2}\Z$, where the sum is finite.

 By $I_i\ast r=0$ for all $i\in\Z$, we can get all the coefficients of the terms containing
$G_j^\pm$ for $j\in\Z+\frac12$ in (\ref{equ-r}) are zero. (In fact
if there exists a non-zero terms $c_{r_1,\cdots
r_n}E_{r_1}\otimes\cdots \otimes E_{r_n}$ in (\ref{equ-r})
containing some $G_{j}^+$ for some $j\in\Z+\frac12$, then we can get
an index sequence $(t_1, t_2, \cdots, t_n)\in (\frac12\Z)^n$ in
which $t_i=r_i$ if $E_{r_i}=G_{r_i}^+$ and zero's in otherwise. We
arrange all such terms (containing some $G_j^+$) by lexicographic
order according to their index sequences and consider the first term
(also called minimal term) $c_{r_1,\cdots r_n}E_{r_1}\otimes
\cdots\otimes E_{r_n}$. Clearly the coefficient of the first term of
all terms in $I_{-1}\ast r$ which including some $G_j^+$ is also
$c_{r_1,r_2,\cdots r_n}$. Then $c_{r_1,r_2,\cdots r_n}=0$. So all
the coefficients of the terms containing $G_j^+$ for
$j\in\Z+\frac12$ are zero. Similar we can consider the terms
including $G_j^-$.)

Similarly, by $I_i \cdot r=0$ for all $i\in\Z$  we can get that all
the coefficients of the terms containing $L_m$ for some $m\in\Z$ are
zero.

Then we can suppose that all terms in (\ref{equ-r}) are only some tensor products of $I_m$ for some $m\in\Z$.
By $G_{q}^+\cdot r=0$ for all $q\in\Z+\frac12$ we can get $r=0$.
This proves the lemma. \QED

As a conclusion of Lemma \ref{lemma2}, one immediately obtains the
following result for the Lie algebra ${\mathcal L}$.
\begin{lemm}\label{colo}
An element $r \in \Im(1 - \tau) \subset {\mathcal L}\otimes
{\mathcal L} $ satisfies CYBE in $(\ref{e-CYBE})$ if and only if it
satisfies MYBE in $(\ref{e-MYBE})$.
\end{lemm}

\section{First cohomology group of $\mathcal L$ in $\mathcal L\otimes \mathcal L$}\setcounter{section}{3}
\setcounter{theo}{0} \setcounter{equation}{0} \vs{6pt}

\begin{theo}\label{theo2}
$\Der(\LL,\VV)=\Inn(\LL,\VV),$ where $\VV=\LL\otimes \LL$,
equivalently, $H^1(\LL,V)=0$.
\end{theo}
\ni{\bf Proof.~} Note that $\VV=\oplus_{i\in\frac{1}{2}\Z}\VV_i$ is
also $\frac{1}{2}\Z$-graded with $\VV_i=\sum_{j+k=i}
\LL_j\otimes\LL_k$, where $i,j,k\in\frac{1}{2}\Z$. We say a
derivation $d\in\Der(\LL ,\VV)$ is {\it homogeneous of degree
$i\in\frac{1}{2}\Z$} if $d(\VV_j) \subset \VV_{i +j}$ for all
$j\in\frac{1}{2}\Z$. Set $\Der(\LL , \VV)_i =\{D\in \Der(\LL , \VV)
\,|\,{\rm deg\,}D =i\}$ for $i\in\frac{1}{2}\Z$.

For any $D\in\Der(\LL,\VV)$, $i\in\frac{1}{2}\Z$, $u\in\LL _j$ with
$j\in\frac{1}{2}\Z$, we can write $D(u)=\sum_{k\in\Z}v_k\in \VV$
with $v_k\in \VV_k$, then we set $D_i(u)=v_{i+j}$. Then
$D_i\in\Der(\LL,\VV)_i$ and
\begin{eqnarray}\label{summable}
\mbox{$D=\sum\limits_{i\in\frac{1}{2}\Z} D_i\
\,\mbox{\ where\ }D_i \in \Der(\LL, \VV)_i,$}
\end{eqnarray}
which holds in the sense that for every $u \in\LL $ only finitely
many $D_i(u)\neq 0,$ and $D(u) = \sum_{i \in\Z} D_i(u)$ (we call
such a sum in (\ref{summable}) {\it summable}).

\begin{clai}
The sum in (\ref{summable}) is finite.
\end{clai}
For any $i\in\Z$, suppose $d_{i}=(v_{i})_{\rm inn}$ for some
$v_{i}\in \VV_{i}$. If $|\{i\,|\,v_i\ne0\}|$ is infinite, then
$d(L_0)=\sum_{i\in\Z}L_0\ast v_{i}=-\sum_{i\in\Z}i v_{i}$ is an
infinite sum, which contradicts $d\in\mbox{Der}(\LL,\VV)$. Thus the
claim and proposition follow. \QED

\begin{clai}\label{clai1}
\rm If \,$i\in\frac{1}{2}\Z\bs\{0\}$, then $D_i\in\Inn(\LL ,\VV)$.
\end{clai}
\ \indent Denote $u=-\frac1i{\ssc\,} D_{i}(L_0)\in \VV_{i}.$ For any
$x_{j}\in \LL_{j},j \in\frac{1}{2}\Z,$  applying $D_{i}$ to
$[L_0,x_{j}]=-j x_{j},$ using $d_{i}(x_j)\in \VV_{i+j}$ and the
action of $L_0$ on $\VV_{i+j}$ is the scalar
$L_0|_{\VV_{i+j}}=-(i+j)$, one has
\begin{equation}\label{equa-add-1}
-(i+j)D_{i}(x_{j}) - (-1)^{[d_{j}][x_{j}]}x_{j}\cdot D_{i}(L_0)=-j
D_{i}(x_{j}),
\end{equation}
i.e., $D_{i}(x_{j})=u_{\rm inn}(x_{j})$, which implies $D_{i}$ is
inner.

From the above we see that Theorem \ref{theo2} follows from the following proposition.

\begin{prop}\label{sub2}
$\Der(\mathcal L, \mathcal L\otimes \mathcal L)_0=\Inn(\mathcal L, \mathcal L\otimes \mathcal L)_0$.
\end{prop}
\noindent{\bf Proof.} For any $D_0\in \Der(\mathcal L, \mathcal
L\otimes \mathcal L)_0$, first we have $D_0(L_0)=0$. In fact, using
(\ref{equa-add-1}) with $i=0$, we obtain $x\ast D_0(L_0)=0,\,\
\forall\,\,x\in\LL_{j},\,j\in\frac{1}{2}\Z$, which together with
Lemma \ref{lemma2} gives $D_0(L_0)=0$.

Now we shall consider $\Der({\mathcal L}_{\bar0}, \mathcal L\otimes \mathcal L)_0$.

Since $D_0({\mathcal L}_{\bar0})\in ({\mathcal L}_{\bar0}\otimes
{\mathcal L}_{\bar0})\oplus (\bigoplus_{X,Y\in \{G^\pm\}} X\otimes
Y)$, where all direct sums are as ${\mathcal L}_{\bar0}$-modules,
from Proposition \ref{prop1}, \ref{der1}, we see that
$$\Der({\mathcal L}_{\bar{0}}, \mathcal L\otimes \mathcal L)_0=\Der({\mathcal L}_{\bar{0}}, {\mathcal L}_{\bar0}\otimes {\mathcal L}_{\bar0})_0.$$

For any $n\in\Z$, one can suppose that
\begin{eqnarray}D_0(L_n)&=&(n\a+\gamma)I_0\otimes
I_n+(n\a^\dag-\gamma)I_n\otimes I_0,\label{dl}\\
D_0(I_n)&=&\b I_0\otimes I_n+\b^\dag I_n\otimes I_0,
\ n\in \Z, \label{di}
\end{eqnarray} for $\a,\a^\dag,\b, \b^\dag, \gamma\in\C$ since $D_0(L_0)=0$.

Since $[L_n, G_r^\pm]=(\frac12n-r)G_{n+r}^\pm, \  [I_n, G_r^\pm]=\pm G_{n+r}^\pm$ and (\ref{dl}),(\ref{di}),
then we can suppose that
$$D_0(G_r^\pm)=\sum_{i\in\Z}c_{r,i}^\pm I_i\otimes G_{r-i}^\pm+\sum_{i\in\Z}d_{r,i}^\pm G_{r-i}^\pm\otimes I_i,$$
where all sums are finite and all coefficients are in $\C$.

Applying $D_0$ to $[L_{2r}, G_r^\pm]=0$, we obtain
\begin{eqnarray}
\!&&\![L_{2r}, D_0(G_r^\pm)]\nonumber\\
\!=\!&&\!-(\pm(2r\a+\gamma)(I_0\otimes G_{3r}^\pm+G_{r}^\pm\otimes I_{2r})
\pm(2r\a-\gamma)(G_{3r}^\pm\otimes I_0+I_{2r}\otimes G_r^\pm)). \label{dg1}
\end{eqnarray}
However, the max$\{i\mid c_{r,i}d_{n,i}\ne 0\}\le0$ and $l$=min$\{i\mid c_{r,i}d_{n,i}\ne 0\}\ge-2r$.
Moreover if $l=-2r$, then $I_{-2r}\otimes G_{5r}^{\pm}$ is in the left side of (\ref{dg1}). It is impossible.
By similar consideration, we obtain $l=0$. Then
we have $D_0(G_r^\pm)=a^\pm I_0\otimes G_r^\pm+b^\pm I_0\otimes G_r^\pm$ for some $a^\pm, b^\pm\in\C$,
and $\a=\a^\dag=\gamma=0$.

Applying $D_0$ to $[G_r^+,G_s^-]=2L_{r+s}+(r-s)I_{r+s}$, we obtain
$\b=\b^\dag=0$ and $a^+=-a^-=b^+=-b^-$. Then replaced $D_0$ by
$D_0+a^+u_{inn}$ for $u=I_0\otimes I_0$, $D_0=0$.\QED

\section{Lie super bialgebra structures on Lie superalgebras}\setcounter{section}{4}
\setcounter{theo}{0} \setcounter{equation}{0} \vs{6pt}

In this section, we shall consider Lie super bialgebra structures on the $N=2$ superconformal Neveu-Schwarz algebra.
First we have
\begin{theo}\label{mainthe}
Every Lie superbialgebra structure on the centerless $N=2$
superconformal Neveu-Schwarz algebra $\LL$ is triangular coboundary.
\end{theo}

First we prove the following lemma.
\begin{lemm} \label{lemma3.4} If $r\in \VV$ satisfies
$x\ast  r\in {\rm Im}(1-\tau)\,(\,\forall\,\,x\in
\LL)$, then $r\in {\rm Im}(1-\tau)$.
\end{lemm}
\ni{\it Proof.~}~Note $\LL \ast  {\rm Im}(1-\tau)\subset
{\rm Im}(1-\tau).$ Write $r=\sum_{i\in\frac{1}{2}\Z}r_i$
with $r_i\in \VV_i$. Obviously, $r\in{\rm Im}(1-\tau)$
if and only if $r_i\in{\rm Im}(1-\tau)$ for all
$i\in\frac{1}{2}\Z.$ Thus without loss of generality, one can
suppose $r=r_i$ is homogeneous.

If $i\in\frac{1}{2}\Z^*$, then $r_i=-\frac1i L_0\ast  r_i\in{\rm
Im}(1-\tau)$. For the case $i=0$, from Lemma 3.5 in \cite{LP} and
$I_0\ast v\in {\rm Im}(1-\tau)$, one can write
\begin{eqnarray*}
r_0\!\!\!&=\mbox{$\sum\limits_{p\in\Z+\frac{1}{2}}$}c_{p}G_{p}^+\otimes
G_{-p}^-+\mbox{$\sum\limits_{p\in\Z+\frac{1}{2}}$}d_{p}G_{p}^-\otimes
G_{-p}^+,
\end{eqnarray*}
where the sum are all finite.  Since the elements of the form
$u_{3,p}:=G_p^+\otimes G_{-p}^--G_{-p}^-\otimes G_{p}^+$
are all in ${\rm
Im}(1-\tau),$ replacing $v$ by $v-u$, where $u$ is a combination of
some $u_{3,r}$,  one can rewritten as
\begin{eqnarray}
r_0\!\!\!=\mbox{$\sum\limits_{p\in\Z+\frac{1}{2}}$}a_{p}G_{p}^+\otimes
G_{-p}^-, \label{wpqr2}
\end{eqnarray}

Since the sum is finite, so there exist $s, t\in\Z+\frac12$, such that
\begin{eqnarray}
r_0\!\!\!=\mbox{$\sum\limits_{s\le p\le t}$}a_{p}G_{p}^+\otimes
G_{-p}^-, \label{wpqr2}
\end{eqnarray} with $a_s, a_t\ne 0$.

However by $I_1\ast r_0\in {\rm Im}(1-\tau)\subset {\rm Ker}(1+\tau)$, we have

\begin{eqnarray}
(1+\tau)\mbox{$\sum\limits_{s\le p\le t}$}a_{p}(G_{p+1}^+\otimes
G_{-p}^--G_{p}^+\otimes
G_{-p+1}^-)=0. \label{wpqr3}
\end{eqnarray}

It is a contradiction with the fact that $a_t\ne 0$. Then one has $a_p=0,\
\forall\,\,p\in\Z+\frac12$. Thus the lemma
follows. \QED
\vskip10pt

\ni{\it Proof of Theorem \ref{mainthe}.} Let $(\LL
,[\cdot,\cdot],\D)$ be a Lie superbialgebra structure on $\LL$. Then
$\D=\D_r$ is defined by (\ref{e-D-r}) for some $r\in\VV_{\bar0} $.
By (\ref{cond2}), ${\rm Im}\, \D\subset{\rm Im}(1-\tau)$. Thus by
Lemma \ref{lemma3.4}, $r\in{\rm Im}(1\otimes 1-\tau)$. Then
(\ref{cond2}), (\ref{e-p2.1}) and Corollary \ref{colo} show that
$c(r)=0$. Thus $(\LL,[\cdot,\cdot],\D)$ is triangular coboundary.
\QED

Now we return to consider Lie superbialgebra structures on
$\hat{\LL}$. As that in Lie algebra case (see Lemma 5.2 in
\cite{NT}), we have following lemma.

\begin{lemm} \label{lemm4}
Let $L$ be a Lie superalgebra such that $L^L$, $(L\otimes L)^L$ and $H^1(L, L)=0$. Then for any
one dimensional central extension $\hat L$ of $L$, there is a linear embedding from
$H^1(\hat L, \hat L\otimes \hat L)$ into $H^1(L, L\otimes L)$. In particular, if $H^1(L, L\otimes L)=0$,
then $H^1(\hat L, \hat L\otimes \hat L)=0$.\QED
\end{lemm}

With Lemma \ref{lemm4} and Theorem \ref{mainthe}, we obtain the
following results.
\begin{theo}\label{mainthe1}
Every Lie superbialgebra structure on the $N=2$ superconformal
Neveu-Schwarz algebra $\hat\LL$ is triangular coboundary.
\end{theo}

\begin{rema}
The above methods can also be used to investigate such structures on
the $N=2$  superconformal twisted, topological, Ramond algebras with
no complicated calculations as in \cite{FL, LFZ,YS}.

\end{rema}

\vskip30pt \centerline{\bf ACKNOWLEDGMENTS}

 \vskip15pt Project is
supported by the NNSF (No. 11071068, 10871057, 11171055), the ZJNSF
(No. D7080080, Y6100148), Qianjiang Excellence Project (No.
2007R10031), the "New Century 151 Talent Project" (2008), the
"Innovation Team Foundation of the Department of Education"  (No.
T200924) of Zhejiang Province, Natural
 Science Foundation of  Jilin province (No.201115006), Scientific
    Research Foundation for Returned Scholars
    Ministry of Education of China and the Fundamental Research Funds for the Central
    Universities.

\end{document}